\DeclareFontFamily{U}{mathb}{\hyphenchar\font45}
\DeclareFontShape{U}{mathb}{m}{n}{ <-6> matha5 <6-7> matha6 <7-8>
mathb7 <8-9> mathb8 <9-10> mathb9 <10-12> mathb10 <12-> mathb12 }{}
\DeclareSymbolFont{mathb}{U}{mathb}{m}{n}

\DeclareMathAccent{\abxring}{0}{mathb}{"38}

\DeclareFontFamily{U}{mathb}{\hyphenchar\font45}
\DeclareFontShape{U}{mathb}{m}{n}{ <-6> matha5 <6-7> matha6 <7-8>
mathb7 <8-9> mathb8 <9-10> mathb9 <10-12> mathb10 <12-> mathb12 }{}
\DeclareSymbolFont{mathb}{U}{mathb}{m}{n}

\documentclass[graybox]{svmult}

\usepackage{type1cm}        
\usepackage{makeidx}         
\usepackage{graphicx}        
\usepackage{multicol}        
\usepackage[bottom]{footmisc}
\usepackage{newtxtext}       
\usepackage{amsfonts}

\usepackage[normalem]{ulem}
\usepackage{soul,xcolor}

\usepackage{caption}

\newcommand\oldtext[1]{}

\makeindex

\usepackage{wrapfig}
\usepackage{graphicx}  
\usepackage{dcolumn}  
\usepackage{bm}       
\usepackage{algorithm,algpseudocode}
\usepackage{enumerate}
 \usepackage{booktabs}
\usepackage{multirow}
\usepackage{makecell}
\usepackage{pgfplots}
\usepackage{tikz}
\usetikzlibrary{patterns}
\usepackage{empheq}
\usepackage{todonotes}
\usepackage{xargs}
\usetikzlibrary{arrows} 
\usepackage{calrsfs}
\usepackage{enumitem}
\usepackage{upgreek}

\definecolor{myblack}{RGB}{53, 53, 53}
\definecolor{myblue}{RGB}{40, 75, 200}
\definecolor{myred}{RGB}{192, 50, 33}
\definecolor{myyellow}{RGB}{255, 166, 48}
\definecolor{mywhite}{RGB}{240, 237, 238}
\definecolor{mygreen}{RGB}{0, 102, 0}
\definecolor{mypurple}{RGB}{150, 0, 180}

\definecolor{green1}{RGB}{9, 82, 86}
\definecolor{green2}{RGB}{8, 127, 140}
\definecolor{green3}{RGB}{6, 167, 125}
\definecolor{green4}{RGB}{79, 109, 122}
\definecolor{green5}{RGB}{192, 214, 223}
\definecolor{violet}{RGB}{26,69,131}

\def \0{\mathbf{0}}			
\def \x{{\mathbf{x}}}			
\def \s{{\mathbf{p}}}			
			
\def \R{{\mathbb{R}}}			
\def \N {\mathbb{N}}

\def \h{{h}}

\def \Rev{{\mathbf{R}}}			
\def \P{{\mathbf{P}}}			
	        
\def \I{{\mathbf{I}}}

\def \x{\vec{x}}				
\def \s{\vec{s}}			
\def \z{\vec{z}}				
\def \c{\vec{c}}

\def \0{\vec{0}}

\def \R{\mathbb{R}}			
\def \N {\mathbb{N}}

\def \Rev{\mathbf{R}}			
\def \P{\mathbf{P}}			
	         
\def \I{\mathbf{I}}

\DeclareMathAlphabet{\mathbfsf}{\encodingdefault}{\sfdefault}{bx}{n}
\DeclareMathAlphabet{\mathsfbfit}{\encodingdefault}{\sfdefault}{bx}{sl}
\DeclareMathAlphabet{\mathbfit}{\encodingdefault}{\rmdefault}{b}{it}
\DeclareMathAlphabet{\mathsfit}{\encodingdefault}{\sfdefault}{m}{sl}

\renewcommand{\vec}[1]{\boldsymbol{#1}}

\usepackage{tikz}
\usetikzlibrary{automata, positioning, arrows, external, matrix}
\usepackage{pgfplots}
\usepgfplotslibrary{fillbetween}
\usepgfplotslibrary{groupplots}
\usepgfplotslibrary{external}

\DeclareMathAlphabet{\pazocal}{OMS}{zplm}{m}{n}
\DeclareMathAlphabet{\bpazocal}{OMS}{cmsy}{b}{n}

\DeclareFontFamily{OT1}{pzc}{}
\DeclareFontShape{OT1}{pzc}{m}{it}{<-> s * [1.10] pzcmi7t}{}
\DeclareMathAlphabet{\mathpzc}{OT1}{pzc}{m}{it}

\def \c{{\mathbf{c}}}	
\usetikzlibrary{matrix}
\usepgfplotslibrary{groupplots, units}
\pgfplotsset{compat=1.3}

\algnewcommand\algorithmiconput{\textbf{Constants:}}
\algnewcommand\algorithmicinput{\textbf{Input:}}
\algnewcommand\algorithmicoutput{\textbf{Output:}}
\algnewcommand{\algorithmicgoto}{\textbf{go to}}

\algnewcommand\Constants{\item[\algorithmiconput]}
\algnewcommand\Input{\item[\algorithmicinput]}\algnewcommand\Output{\item[\algorithmicoutput]}\algnewcommand{\Goto}[1]{\algorithmicgoto~\ref{#1}}

\algrenewcommand\ALG@beginalgorithmic{\footnotesize}
\algrenewcommand\alglinenumber[1]{\scriptsize #1:}

\usetikzlibrary{decorations.pathmorphing,shapes}

 \pgfplotsset{width=7cm,compat=1.3}

\def\blfootnote{\gdef\@thefnmark{}\@footnotetext}

\begin{document}

\title*{On the use of hybrid coarse-level models in multilevel minimization methods}
\author{Alena Kopani\v{c}\'akov\'a}
\institute{Alena Kopani\v{c}\'akov\'a 	\at Brown University, USA,  \email{alena.kopanicakova@brown.edu}
}
\maketitle

\abstract{Solving large-scale nonlinear minimization problems is computationally demanding. Nonlinear multilevel minimization (NMM) methods explore the structure of the underlying minimization problem to solve such problems in a computationally efficient and scalable manner. The efficiency of the NMM methods relies on the quality of the coarse-level models. Traditionally, coarse-level models are constructed using the additive approach, where the so-called~$\tau$-correction enforces a local coherence between the fine-level and coarse-level objective functions. In this work, we extend this methodology and discuss how to enforce local coherence between the objective functions using a multiplicative approach. Moreover, we also present a hybrid approach, which takes advantage of both, additive and multiplicative, approaches. Using numerical experiments from the field of deep learning, we show that employing a hybrid approach can greatly improve the convergence speed of NMM methods and therefore it provides an attractive alternative to the almost universally used additive approach. }

\blfootnote{\noindent \textbf{Funding:} This work was supported by the Swiss National Science Foundation under the project "Multilevel training of DeepONets --- multiscale and multiphysics applications" (grant no.~206745). 
}

\section{Introduction}
\label{sec:intro}
We consider the following minimization problem:
\begin{align}
\min_{\x \in \R^n} f(\x),
\label{eq:min_problem}
\end{align}
where $f:\R^n \rightarrow \R$ is a bounded, twice continuously differentiable objective function and $n \in \N$ is typically very large. 
Our goal is to minimize~\eqref{eq:min_problem} using a nonlinear multilevel minimization (NMM) method, e.g.,~MG-OPT~\cite{Nash2000multigrid} or RMTR~\cite{Gratton2008recursive}.  
The main idea behind NMM methods is to employ a hierarchy of so-called coarse-level objective functions, denoted by $\{f^{\ell} \}_{\ell=1}^L$, where $L>1$.
These functions are typically obtained by exploring the structure of the underlying minimization problem, e.g., by discretizing the underlying infinite-dimensional problem with a varying discretization parameter. 
During the solution process, the functions~$\{f^{\ell} \}_{\ell=1}^L$ are utilized in order to construct the search-directions for the minimization problem at hand in a computationally efficient manner.

The overall efficiency of NMM methods relies on the ability of the coarse-level objective functions~$\{f^{\ell} \}_{\ell=1}^L$ to approximate the function $f$ well.  
Indeed, the convergence theory of the majority of NMM methods requires that the local behavior of the coarse-level objective functions is at least first-order coherent with the local behavior of~$f$. 
The coherence is commonly ensured by employing the so-called $\tau$-correction~\cite{brandt1977multi}, which corrects the coarse-level objective function $f^{\ell}$ in an additive manner. 
Although this approach is almost universally employed in the multilevel literature, other approaches were also considered, e.g.,~a second-order additive correction approach~\cite{Yavneh2006, Gratton2008recursive,kopanivcakova2020recursive}, or Galerkin-based coarse-level models~\cite{ho2019newton, Gratton2008recursive}. 
In this work, we explore techniques from the surrogate-based/multi-fidelity optimization~\cite{forrester2009recent} in order to construct the first-order coherent coarse-level models in the context of NMM methods.
In particular, we discuss how to correct functions~$\{f^{\ell} \}_{\ell=1}^L$ using additive, multiplicative, and hybrid approaches.

\section{Nonlinear multilevel minimization framework}
In this work, we minimize~\eqref{eq:min_problem} using the NMM method.
To this aim, we consider a hierarchy of $L$ levels. 
Each level~$\ell = 1, \dots, L$ is associated with some model~${h^{\ell}: \R^{n^{\ell}} \rightarrow \R}$, where we assume that~$h^{\ell-1}$ is computationally cheaper to minimize than~$h^{\ell}$ and that $n^{\ell-1} < n^{\ell}$. 
As we will discuss in Section~\ref{sec:models}, the models~$\{ h^{\ell} \}_{\ell=1}^{L}$ are constructed during the minimization process by correcting the objective functions~$\{ f^{\ell} \}_{\ell=1}^{L}$ by taking into account the knowledge of the current iterate.
Through this work, we assume that~$h^L:=f^L:=f$.
Transfer of the data between different levels of the multilevel hierarchy is performed using the prolongation operator~${\I^{\ell+1}_{\ell}: \R^{n^{\ell}} \rightarrow \R^{n^{\ell+1}}}$, and the restriction operator~${\Rev_{\ell+1}^{\ell}: \R^{n^{\ell+1}} \rightarrow \R^{n^{\ell}}}$, where~$\Rev_{\ell+1}^{\ell} = (\I^{\ell+1}_{\ell})^T$. 
Moreover, we also employ the projection operator~${\P_{\ell+1}^{\ell}: \R^{n^{\ell+1}} \rightarrow \R^{n^{\ell}}}$ to transfer iterates from the level~$\ell+1$ to~$\ell$.
The operator~$\P_{\ell+1}^{\ell}$ is constructed such that $\x^{\ell} = \P_{\ell+1}^{\ell}(\I_{\ell}^{\ell+1}\x^{\ell})$, for any~$\x^{\ell} \in \R^{n^{\ell}}$.

Using the aforementioned definitions, we now describe a generic NMM method in the form of a V-cycle, summarized in Algorithm~\ref{alg:nonlinear_method}. 
During the description, we use a superscript to denote the level and a subscript to denote the iteration index.
Starting from the finest level, $\ell= L$, and initial guess $\x_0^{\ell}$, the NMM method performs $\mu_s$ nonlinear smoothing steps to approximately minimize model~$h^{\ell}$. 
The choice of the nonlinear smoother depends on the particular choice of NMM method.
For instance, one can employ a first-order method equipped with a line-search or trust-region globalization strategy if a variant of multilevel line-search or trust-region method is considered. 
The outcome of this minimization process, iterate $\x_{\mu_s}^{\ell}$, is then used to construct a coarse-level model $h^{\ell-1}$ and initial guess~$\x^{\ell-1}_{0} = \P_{\ell}^{\ell-1} \x^{\ell}_{\mu_s}$.
This process is repeated recursively until the coarsest level is reached.

On the coarsest level, $\ell=1$, an NMM method approximately minimizes~$h^{\ell}$ using $\mu_c$ steps of a nonlinear solution strategy, giving rise to~$\x^{\ell}_{*}$.
Afterwards, the prolongated coarse-level correction~$\s^{\ell+1}_{\mu_s+1}:= {\I_{\ell}^{\ell+1} (\x^{\ell}_{*} - \x^{\ell}_0)}$ is used to update the current iterate~$\x^{\ell+1}_{\mu_s}$ on level~$\ell+1$. 
However, before this update is performed, the correction~$\s^{\ell}_{\mu_s+1}$ has to undergo some convergence control.
The type of convergence control again depends on the particular type of the NMM method. 
For example, if the multilevel trust-region method is used, then~$\s^{\ell+1}_{\mu_s+1}$ is required to provide a decrease in $h^{\ell+1}$ to be accepted by the algorithm. 
If a variant of a line-search method is used, then an appropriate step size has to be determined. 
In the end, the algorithm performs $\mu_s$ post-smoothing steps, starting from $\x^{\ell+1}_{\mu_s+1}$ and giving rise to $\x^{\ell+1}_{*}$. 
This process is again repeated on all levels until the finest level is reached.

\begin{algorithm}
\footnotesize
\caption{NMM($\ell, \ h^{\ell}, \  \x_0^{\ell}$)}
\label{alg:nonlinear_method}
\begin{algorithmic}[1]
\Require{$\ell \in \N, h^{\ell}:\R^{n^{\ell}} \rightarrow \R, \  \x_0^{\ell} \in \R^{n^{\ell}}$ and $\mu_s, \mu_c \in \N$}
\State  $\x_{\mu_s}^{\ell}$ = Nonlinear\_smoothing($h^{\ell}, \ \x_0^{\ell}, \ \mu_s$)
\State  Construct $h^{\ell-1}$ using $\x_{\mu_s}^{\ell}$, and $\nabla h^{\ell}(\x_{\mu_s}^{\ell})$
\If{$\ell=2$}
\State
$\x_{*}^{\ell-1}$ = Nonlinear\_solve($h^{\ell-1}, \ \P_{\ell}^{\ell-1} \bm{\x}_{\mu_s}^{\ell}, \ \mu_c$)
\Else
\State
$\x_{*}^{\ell-1}$ = NMM($\ell-1,  h^{\ell-1}, \  \P_{\ell}^{\ell-1} \bm{\x}_{\mu_s}^{\ell}$)
\EndIf
\State $\x_{\mu_s+1}^{\ell}$ =  Convergence\_control($ h^{\ell}, \x_{\mu_s}^{\ell},  \bm{I}_{\ell-1}^{\ell} (\x_{*}^{\ell-1} -  \P_{\ell}^{\ell-1} \x_{\mu_s}^{\ell}))$
\State  $\x_{*}^{\ell}$ = Nonlinear\_smoothing($h^{\ell}, \ \x_{\mu_s+1}^{\ell}, \ \mu_s$)
\State \Return $\x_{*}^{\ell}$
\end{algorithmic}
\end{algorithm}

\section{Construction of coarse-level models}
\label{sec:models}
On each level~$\ell$, the NMM methods minimize the model~${\h^{\ell}:\R^{n^{\ell}} \rightarrow \R}$ approximately. 
The result of this minimization, the iterate~$\x_{*}^{\ell}$, is then used to construct the search direction for the minimization on the next finer level. 
As a consequence, the overall efficiency of NMM methods depends on the capabilities of the models~$\{\h^{\ell}\}_{\ell=1}^L$ to approximate $f$ as accurately as possible.

Given an initial guess~${\x^{\ell}_0=\P_{\ell+1}^{\ell} \x^{\ell+1}_{\mu_s}}$, the model~$\h^{\ell}$ is constructed during each V-cycle by correcting the  function$f^{\ell}$, such that 
the following condition holds: 
\begin{align}
\nabla \h^{\ell}(\x^{\ell}_{0}) = \Rev^{\ell}_{\ell+1} \nabla \h^{\ell+1} (\x^{\ell+1}_{\mu_s}).
\label{eq:first_order_consistency}
\end{align} 
This ensures that~$\h^{\ell}$ and~$\h^{\ell+1}$ are locally first-order coherent and that the following relation holds: ${\langle \nabla \h^{\ell}(\x^{\ell}_{0}) , \s^{\ell} \rangle = \langle   \nabla \h^{\ell+1} (\x^{\ell+1}_{\mu_s}) , \I_{\ell}^{\ell+1} \s^{\ell} \rangle}$.
In this work, we discuss three different approaches for constructing models $\{ h^{\ell}\}_{\ell=1}^{L}$, namely additive, multiplicative and hybrid. 
Our discussion considers only the first-order coherent models, constructed using the Taylor approximation of the associated correction function. 
However, models enforcing higher-order coherency as well as different approximations of the correction function could also be considered.

\subsection{An additive approach}
\label{sec:additive_approach}
Using the additive approach, the coarse-level model~${\h_{\text{add}}^{\ell}: \R^{n^{\ell}} \rightarrow \R}$ is obtained by correcting the low-cost function~$f^{\ell}$ as follows
\begin{align}
\h^{\ell}_{\text{add}}(\x^{\ell}) = f^{\ell}(\x^{\ell}) + \gamma_{\text{add}}^{\ell}(\x^{\ell}),
\label{eq:model_additive_generic}
\end{align}
where the additive correction function ${\gamma_{\text{add}}^{\ell}: \R^{n^{\ell}} \rightarrow \R}$ accounts for the difference between  the value of~$f^{\ell}$ and the fine-level model~$\h^{\ell+1}$, i.e., 
\begin{align}
\gamma_{\text{add}}^{\ell}(\x^{\ell}) := \h^{\ell+1}( \x^{\ell+1}_{\mu_{s}}) -  f^{\ell}(\x^{\ell}).
\label{eq:gamma_corr}
\end{align}

Unfortunately, the evaluation of~$\gamma_{\text{add}}^{\ell}$ at any given $\x^{\ell}$ requires an evaluation of the fine-level model~$\h^{\ell+1}$ at $\I_{\ell}^{\ell+1} \x^{\ell}$.
As a consequence, numerical computations involving~$h^{\ell}_{\text{add}}$ are computationally more demanding than computations performed using~$h^{\ell+1}$ directly.
To ease the computational burden, we evaluate~$\gamma_{\text{add}}^{\ell}$ exactly only at the initial coarse-level iterate~${\x^{\ell}_0=\P^{\ell+1}_{\ell} \x^{\ell+1}_{\mu_s}}$.
Thus, we impose 
\begin{align*}
\gamma_{\text{add}}^{\ell}(\x^{\ell}_0) := \h^{\ell+1}(\x^{\ell+1}_{\mu_s}) - f^{\ell}(\x^{\ell}_0),
\end{align*}
only at~$\x^{\ell}_0$. 
For any other~$\x^{\ell}$, we approximate the correction function~$\gamma_{\text{add}}^{\ell}$ by means of the first-order Taylor approximation, defined around~$\x^{\ell}_0$ as follows
\begin{align*}
\tilde{\gamma}_{\text{add}}^{\ell}(\x^{\ell}) = \gamma_{\text{add}}^{\ell}(\x^{\ell}_0) + \langle \nabla \gamma_{\text{add}}^{\ell}(\x^{\ell}_0),  \ \x^{\ell} -\x^{\ell}_0 \rangle.
\end{align*}
Replacing~$\gamma^{\ell}_{\text{add}}$ with~$\tilde{\gamma}^{\ell}_{\text{add}}$ in~\eqref{eq:model_additive_generic} gives rise to
\begin{align}
\h^{\ell}_{\text{add}}(\x^{\ell}) := f^{\ell}(\x^{\ell}) + 
 \h^{\ell+1}(\x^{\ell+1}_{\mu_s}) - f^{\ell}(\x^{\ell}_0) +
 \langle 
 \nabla \gamma_{\text{add}}^{\ell}(\x^{\ell}_0), \x^{\ell} -\x^{\ell}_0 \rangle,
\label{eq:coarse_objective_first_order}
\end{align}
where 
\begin{align}
 \nabla \gamma_{\text{add}}^{\ell}(\x^{\ell}_0) := 
\Rev^{\ell}_{\ell +1} \nabla h^{\ell+1} (\x^{\ell+1}_{\mu_s})
- \nabla f^{\ell}(\x^{\ell}_{0}). 
\label{eq:delta_g}
\end{align}
Note, the quantity~$ \h^{\ell+1}(\x^{\ell+1}_{\mu_s}) - f^{\ell}(\x^{\ell}_0)$
enforces zeroth-order coherence between~$\h^{\ell+1}$ and~$\h_{\text{add}}^{\ell}$ at~$\x^{\ell+1}_{\mu_s}$ and~$\x^{\ell}_0$, respectively, i.e.,~${\h_{\text{add}}^{\ell}(\x^{\ell}_0) = \h^{\ell+1}(\x^{\ell+1}_{\mu_s})}$. 
However, this term does not affect the evaluation of the derivatives of~$\h^{\ell}_{\text{add}}$, and therefore it is often neglected in practice.
We also point out that the term~$ \nabla \gamma_{\text{add}}^{\ell}(\x^{\ell}_0)$, known in the multilevel literature as $\tau$-correction, ensures that condition~\eqref{eq:first_order_consistency} holds.

\subsection{A multiplicative approach}
\label{sec:multiplicative}
Optimization methods that exploit multiple fidelities often employ multiplicative correction functions~\cite{forrester2009recent}. 
In this case, the low-cost approximation~$f^{\ell}$ associated with level~$\ell$ is made coherent with the model~$h^{\ell+1}$ as follows:
\begin{align}
\h^{\ell}_{\text{mult}}(\x^{\ell}) = \gamma^{\ell}_{\text{mult}}(\x^{\ell}) f^{\ell}(\x^{\ell}),
\label{eq:mult_gg}
\end{align}
where ${\gamma^{\ell}_{\text{mult}}:\R^{n^{\ell}} \rightarrow \R}$ is a multiplicative correction function, given as
\begin{align}
\gamma^{\ell}_{\text{mult}}(\x^{\ell}) := \frac{\h^{\ell+1}( \I_{\ell}^{\ell+1}\x^{\ell})+ \kappa}{f^{\ell}(\x^{\ell})+ \kappa},
\label{eq:mul_generic}
\end{align}
where $\kappa \approx \epsilon$ ensures numerical stability as the value of~$f^{\ell}(\x^{\ell})$ approaches zero.

Similar to the additive approach, evaluating~$\gamma^{\ell}_{\text{mult}}$ precisely at all coarse-level iterates is computationally expensive. 
Therefore, we impose~\eqref{eq:mul_generic} only at~$\x^{\ell}_0$, i.e.,
\begin{align*}
\gamma^{\ell}_{\text{mult}}(\x^{\ell}_0) := \frac{\h^{\ell+1}(\x^{\ell+1}_{\mu_s})+\kappa}{f^{\ell}(\x^{\ell}_0)+\kappa},
\end{align*}
where we explored that~$\x^{\ell+1}_{\mu_s} = \I_{\ell}^{\ell +1} \x^{\ell}_0$.
At any other iterate~$\x^{\ell}$,  we approximate~$\gamma^{\ell}_{\text{mult}}$ by means of the first-order Taylor approximation, defined around~$\x_0^{\ell}$ as
\begin{align}
\tilde{\gamma}^{\ell}_{\text{mult}}(\x^{\ell}) = \gamma^{\ell}_{\text{mult}}(\x^{\ell}_0) + \langle \nabla \gamma^{\ell}_{\text{mult}}(\x^{\ell}_0),  \ \x^{\ell} - \x^{\ell}_0 \rangle.
\label{eq:first_app_betta_mult}
\end{align}
Replacing  $\gamma^{\ell}_{\text{mult}}$ with~$\tilde{\gamma}^{\ell}_{\text{mult}}$ in~\eqref{eq:mult_gg} then gives rise to the first-order coherent model
\begin{align}
\h^{\ell}_{\text{mult}}(\x^{\ell}) := \tilde{\gamma}^{\ell}_{\text{mult}}(\x^{\ell}) \ \  f^{\ell}(\x^{\ell}).
\label{eq:mult_first}
\end{align}
The numerical evaluation of~$\tilde{\gamma}^{\ell}_{\text{mult}}$ amounts to
\begin{equation*}
\begin{aligned}
\tilde{\gamma}^{\ell}_{\text{mult}}(\x^{\ell}) := & \frac{\h^{\ell+1}(\x^{\ell+1}_{\mu_s})+\kappa}{f^{\ell}(\x^{\ell}_0)+ \kappa}  +  \langle \nabla \gamma^{\ell}_{\text{mult}} (\x_0^{\ell}), \ \x^{\ell} - \x^{\ell}_0\rangle,  
\end{aligned}
\end{equation*}
where $\nabla \gamma^{\ell}_{\text{mult}}(\x_0^{\ell}) $ is given by
\begin{align*}
\nabla \gamma^{\ell}_{\text{mult}}(\x_0^{\ell})  := \frac{1}{f^{\ell}(\x^{\ell}_0)} \big(\Rev^{\ell}_{\ell+1} \nabla \h^{\ell+1}({\x^{\ell+1}_{\mu_s}}) \big)  - \frac{\h^{\ell+1}(\x^{\ell+1}_{\mu_s})}{(f^{\ell}(\x^{\ell}_0))^2} \nabla f^{\ell}(\x^{\ell}_0).
\end{align*}
Straightforward calculations show that model~$\h^{\ell}_{\text{mult}}$, defined by~\eqref{eq:mult_first}, is zeroth-order and first-order coherent with~$h^{\ell+1}$ at~$\x^{\ell}_0$ and~$\x^{\ell+1}_{\mu_s}$, respectively. 

\subsection{A hybrid approach}
\label{sec:hybrid_approach}
\begin{figure}[t]
\begin{tikzpicture}
\begin{groupplot}[group style={group size= 2 by 1},
xmode=normal, 
ymode=normal,  
width=0.55\textwidth, 
height=0.3\textwidth]

\nextgroupplot[ylabel={}, 
xlabel={}, 
ymajorgrids=true, 
xmajorgrids=true, 
tick label style={font=\scriptsize}, 
ymax=20, 
ymin=-10,
xmin=0, 
xmax=10, ]
    \addplot[color = myblue, very thick] table [x=x, y=fine, col sep=comma] {functions_point2.csv};    
    \addplot[color = myred,  very thick] table [x=x, y=coarse, col sep=comma] {functions_point2.csv};    
  
     \addplot[color = mypurple,  very thick, dotted] table [x=x, y=additive_first, col sep=comma] {functions_point1.csv};    
    \addplot[color = green3, very thick, dashed] table [x=x, y=mult_first, col sep=comma] {functions_point1.csv};        
    \addplot[color = myyellow,  very thick, dashdotted] table [x=x, y=mixed, col sep=comma] {functions_point1.csv};        
    
  \addplot[very thick, only marks, mark=x, mark size=4pt] table [x=x_point1,   y=y_point1, col sep=comma] {functions_point1.csv};       
  
\node at (axis cs: 3.0, 3.1) {\footnotesize $x^{L-1}_{0}$};  

\nextgroupplot[ylabel={}, 
xlabel={}, 
ymajorgrids=true, 
xmajorgrids=true, 
tick label style={font=\scriptsize}, 
ymax=10, 
ymin=-8,
xmin=0, 
xmax=10, ]
   
    \addplot[color = myblue, very thick] table [x=x, y=fine, col sep=comma] {functions_point2.csv};   
    \label{pgfplots:Plot1}
     
    \addplot[color = myred,  very thick] table [x=x, y=coarse, col sep=comma] {functions_point2.csv};    
    \label{pgfplots:Plot2}    
    
    \addplot[color = mypurple,  very thick, dotted] table [x=x, y=additive_first, col sep=comma] {functions_point2.csv};    
    \label{pgfplots:Plot3}        
    
    \addplot[color = green3, very thick, dashed] table [x=x, y=mult_first, col sep=comma] {functions_point2.csv};
    \label{pgfplots:Plot4}        
                    
    \addplot[color = myyellow,  very thick, dashdotted] table [x=x, y=mixed, col sep=comma] {functions_point2.csv};
    \label{pgfplots:Plot5}                 
   
\addplot[very thick, only marks, mark=x, mark size=4pt] table [x=x_point2,   y=y_point2, col sep=comma] {functions_point1.csv};    
\node at (axis cs: 6.25, -2.4) {\footnotesize $x^{L-1}_{0}$};
\end{groupplot}
\matrix[
    matrix of nodes,
    anchor=north west,
    fill=white,draw,
    inner sep=0.3em,
    node font=\scriptsize,
  ]
  at (0.0, 2.7){
    \ref{pgfplots:Plot1}& $f^{L}$ &  \hspace{0.38cm} & \ref{pgfplots:Plot2}& $f^{L-1}$ &  \hspace{0.38cm} &
    \ref{pgfplots:Plot3}& $h^{L-1}_{\text{add}}$ &  \hspace{0.38cm} & \ref{pgfplots:Plot4}& $h^{L-1}_{\text{mult}}$ &   \hspace{0.38cm} &
    \ref{pgfplots:Plot5}& $h^{L-1}_{\text{mix}} (\text{MFV})$ \\};
\end{tikzpicture}
\caption{
Coarse-level models constructed around~$x^{L-1}_0=2.5$ and~$x^{L-1}_0=6.0$.
}
\label{fig:obj_example_ex1}
\end{figure}
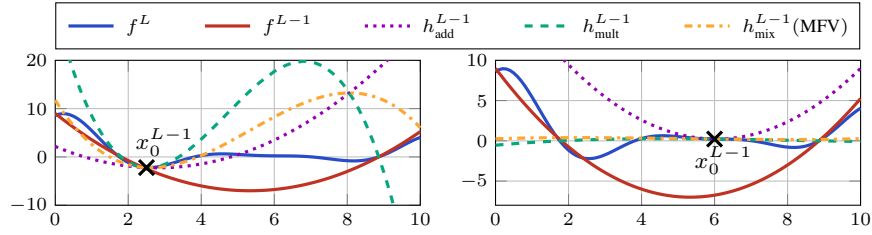

From a computational point of view, additive and multiplicative approaches are comparable. 
However, their behavior is very different. 
The additive approach adds new terms to~$f^{\ell}$, which can be interpreted as uniform translation (zeroth-order), and rotation (first-order) of the function graph; see also Fig.~\ref{fig:obj_example_ex1}. 
In contrast, the multiplicative approach introduces skewing, which might not be desirable if~$f$ and $f^{\ell}$ are in good agreement, at least locally. 
However, if functions $f^{\ell}$ and $f$ are not in good agreement, then additional skewing can be beneficial~\cite{fischer2018bayesian}, e.g.,~if the polynomial order of $f$ is higher than the polynomial order of~$f^{\ell}$.
Moreover, multiplication of $f^{\ell}$ with~$\tilde{\gamma}^{\ell}_{\text{mult}}$ can introduce new minima. 
For instance, let us suppose that~$f^{\ell}$ is a second-order polynomial. 
Its multiplication with~$\tilde{\gamma}^{\ell}_{\text{mult}}$ increases the order of the polynomial, i.e.,~we obtain a model~$h^{\ell}_{\text{mult}}$ which is quartic and has, in general, more minima than quadratic function.

In general, it is not known a priori whether the additive or the multiplicative model is more suitable for a given optimization problem.
To overcome this difficulty, a hybrid approach~\cite{gano2005hybrid} can be employed.
A coarse-level model~$h^{\ell}_{\text{mix}}$ is then obtained as a convex combination of the additive~$\h^{\ell}_{\text{add}}$ and the multiplicative~$\h^{\ell}_{\text{mult}}$ models, i.e., 
\begin{align}
\h^{\ell}_{\text{mix}}(\x^{\ell}) := w^{\ell}_{\text{add}} \ \h^{\ell}_{\text{add}}(\x^{\ell}) + w^{\ell}_{\text{mult}} \ \h^{\ell}_{\text{mult}}(\x^{\ell}),
\label{eq:hybrid_function}
\end{align}
where $w^{\ell}_{\text{add/mult}} \in \R$ and $w^{\ell}_{\text{add}}+w^{\ell}_{\text{mult}}=1$. 
In order to maximize the approximation properties of~$\h^{\ell}_{\text{mix}}$, the weights $w^{\ell}_{\text{add}}, w^{\ell}_{\text{mult}}$ have to be chosen carefully. 
Below, we describe two different strategies for selecting the values~$w^{\ell}_{\text{add}}$ and $w^{\ell}_{\text{mult}}$. 

\subsubsection{Matching function values (MFV) at the previously evaluated fine-level iterate}
Following~\cite{eldred2006formulations}, the weights $w^{\ell}_{\text{add}}, w^{\ell}_{\text{mult}}$ can be selected by matching the function value at the previously evaluated fine-level iterate, denoted by~$\x^{\ell+1}_p$, as in
\begin{align}
w^{\ell}_{\text{add}} = \frac{h^{\ell+1}(\x^{\ell+1}_p) - \h^{\ell}_{\text{mult}}(\x^{\ell}_0)  }{ \h^{\ell}_{\text{add}}(\x^{\ell}_0) - \h^{\ell}_{\text{mult}}( \x^{\ell}_0) } \qquad \quad \text{and} \qquad \quad w^{\ell}_{\text{mult}} = 1 - w^{\ell}_{\text{add}}.
\label{eq:w_form}
\end{align}
From a computational point of view, evaluating~\eqref{eq:w_form} is cheap as~$h^{\ell+1}(\x^{\ell+1}_p)$ is readily available, for instance from the 
$\mu_{s}-1$ pre-smoothing step performed on level $\ell+1$.

\subsubsection{Bayesian updating approach}
To maximize the approximation properties of~$h^{\ell}_{\text{mix}}$, it might be beneficial to take into account the history of the $d^{\ell}$ previously evaluated fine-level iterates~\cite{fischer2018bayesian}. 
Therefore, we consider dataset~$\pazocal{D}^{\ell} = \{(h^{\ell+1}(\x_p^{\ell+1}), h^{\ell}_{\text{add}}(\P_{\ell+1}^{\ell} \x_{p}^{\ell+1} ), h^{\ell}_{\text{mult}}(\P_{\ell+1}^{\ell} \x_{p}^{\ell+1} ) \}_{p=1}^{d^{\ell}}$, where each sample contains the function value of~$\h^{\ell+1}$ at~$\x_p^{\ell+1}$, as well as the function values of the coarse-level models $h^{\ell}_{\text{add/mult}}$ obtained at~$\P_{\ell+1}^{\ell} \x_{p}^{\ell+1}$. 
In this work, we construct~$\pazocal{D}^{\ell}$ by taking into account the last $d^{\ell}$ iterates which were transferred from level~$\ell+1$ to level~$\ell$.
For example, if~$d^{\ell}=3$, then $\pazocal{D}^{\ell}$ is constructed by taking into account the iterate~$\x_p^{\ell+1} = \x_{\mu_s}^{\ell+1}$, obtained as a result of the pre-smoothing step during the previous three V-cycles. 
For simplicity, we use the notation~$d^{\ell}=\infty$ to denote all previous V-cycles.

Having constructed the dataset~$\pazocal{D}^{\ell}$, we can now employ the Bayesian posterior updating approach~\cite{fischer2018bayesian} to determine the values of~$w^{\ell}_{\text{add}/\text{mult}}$. 
Starting from~$w^{\ell}_{\text{add}/\text{mult}}=0.5$, the weights are updated every time the model~$h^{\ell}$ is constructed as follows:
\begin{align}
w^{\ell}_{\text{add}/\text{mult}} = \frac{w^{\ell}_{\text{add}/\text{mult}}  \psi^{\ell}_{\text{add}/\text{mult}} }{w^{\ell}_{\text{mult}/\text{add}}  \psi^{\ell}_{\text{mult}/\text{add}}     + w^{\ell}_{\text{add}/\text{mult}}   \psi^{\ell}_{\text{add}/\text{mult}} }.
\label{eq:weights}
\end{align}
The model likelihoods~$\psi^{\ell}_{\text{add}/\text{mult}}$ in~\eqref{eq:weights} are evaluated as 
\begin{align}
 \psi^{\ell}_{\text{add}/\text{mult}}(\x^{\ell}) = \big(2 \pi {\sigma}^2_{\text{add}/\text{mult}} \big)^{-{d^{\ell}}/2} \exp(-{d^{\ell}}/2 ),
\end{align}
and the maximum likelihood estimator of the model variance is given by
\begin{align}
{\sigma}^2_{\text{add}/\text{mult}} = \frac{1}{d^{\ell}} \sum_{p=1}^{d^{\ell}}	(h^{\ell+1}(\x_p^{\ell+1}) -  h^{\ell}_{\text{add}/\text{mult}}(\P^{\ell}_{\ell+1} \x_{p}^{\ell+1} )  ).
\end{align}

\section{Numerical results and discussion}
\label{sec:num_results}
In this section, we investigate the influence of different coarse-level models on the performance of the NMM method using numerical examples from the field of supervised learning, namely classification using ResNets~\cite{he2016deep}.
Given a dataset~${\pazocal{S}=\{(\z_s, \c_s )\}_{s=1}^{n_s}}$, where $\z_s \in \R^{n_{in}}$ and $\c_s \in \R^{n_{out}}$, our goal is to find parameters $\x \in \R^{n}$ of a ResNet, defined as~${\text{RN}:\R^{n_{in}} \times \R^{n}\rightarrow \R^{n_{out}}}$, by solving the following minimization problem:
\begin{align}
\min_{\x \in \R^{n}} f(\x) := \frac{1}{n_s} \sum_{s=1}^{n_s} g(\text{RN}(\z_s, \x), \c_s),
\label{eq:min_problem_ex}
\end{align}
where $g$ denotes the cross-entropy loss function.

\begin{figure}[t]
\centering
\caption{\emph{Left:} Blobs, Smiley, and Spiral datasets (\emph{Top} to \emph{Down}). 
Each class is illustrated by different color.
\emph{Right:} The average computational cost of the RMTR method ($4$ levels) for training ResNets.
Averages are obtained from~$5$ independent runs.}
\label{fig:results}
\hspace{0.3cm}
\begin{minipage}{0.2\linewidth}
\includegraphics[scale=0.1]{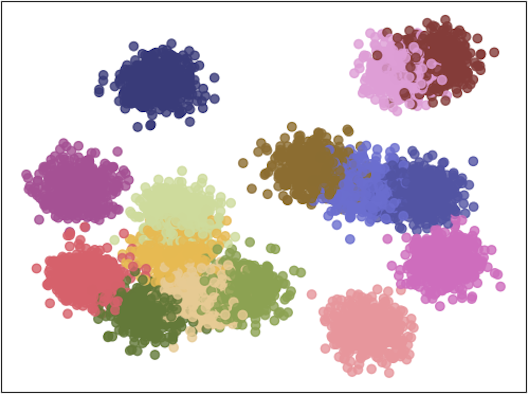}\\
\includegraphics[scale=0.1]{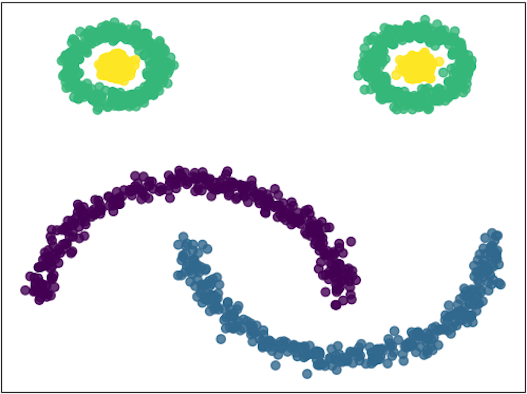}\\
\includegraphics[scale=0.1]{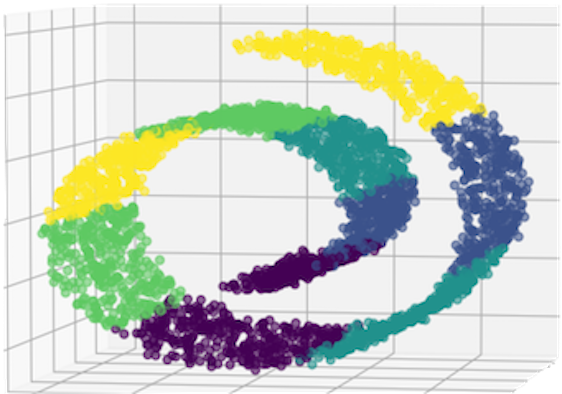}
\end{minipage}
\hspace{0.3cm}
\begin{minipage}{0.7\linewidth}
\small
  \begin{tabular}{ l|c|c|c}
    \toprule
    {Model/Example}                 		& Blobs 			&   Smiley 		&   Spiral       		      \\ \hline \hline

$h_{\text{add}}$			& $29\pm5.3\%$  & $676 \pm11.2\%$ & $203 \pm12.3\%$  \\
$h_{\text{mult}}$			& $32\pm6.1\%$  & $85 \pm15.1\%$ & $153 \pm15.9\%$   \\ \hline
$h_{\text{mix}} (w=0.5) $		& $38 \pm4.8\%$  & $404 \pm10.3\%$ & $297 \pm11.3\%$  \\ \hline
$h_{\text{mix}} (\text{MFV}) $	& $25\pm4.2\%$  & $352 \pm6.5\%$ & $\bm{123 \pm7.1\%}$  \\ \hline
$h_{\text{mix}} (d^{\ell} = 5)$		& $25\pm3.4\%$  & $514 \pm6.3\%$ & $197 \pm6.8\%$  \\
$h_{\text{mix}} (d^{\ell} = 20)$		& $\bm{24\pm2.9\%}$  & $471 \pm7.7\%$ & $156 \pm7.4\%$  \\
$h_{\text{mix}} (d^{\ell} = \infty)$	& $25\pm3.8\%$  & $\bm{301 \pm6.9\%}$ & $126 \pm9.9\%$  \\

\end{tabular}
\end{minipage}
\end{figure}

Since~\eqref{eq:min_problem_ex} is a non-convex function, we choose the NMM method to be a variant of the RMTR method~\cite{Gratton2008recursive}. 
The multilevel hierarchy and transfer operators are constructed by leveraging the fact that the ResNet can be interpreted as a forward Euler discretization of an ordinary differential equation; see~\cite{kopanicakova2022globally, gaedke2021multilevel} for details. 
Here, we construct a hierarchy of ResNets by uniformly refining a ResNet with three layers three times.
Fig.~\ref{fig:results} demonstrates the computational cost of the RMTR method with respect to different coarse-level models for three different datasets. 
As we can observe, the choice of the coarse-level model has a significant impact on the overall efficiency of the mutlilevel method. 
For all three examples, hybrid approaches outperform purely additive and multiplicative ones.
In terms of hybrid models, we observe that the Bayesian approach performs similar, or superior to MFV, especially if all prior fine-level iterates are considered ($d^{\ell}=\infty$). 

Given our (limited) numerical experience, we believe that employing hybrid, and possibly other types of novel coarse-level models,  provides a promising future direction for improving the efficiency and the reliability of NMM methods.

\bibliographystyle{spmpsci}
\bibliography{biblio.bib}

\end{document}